\documentclass[1pt,notitlepage,twoside,a4paper]{amsart}

\usepackage{amsmath,amssymb,enumerate}

\usepackage{epsfig,fancyhdr,color}

\usepackage{amssymb}
\usepackage{amsmath,amsthm}    
\usepackage{latexsym}
\usepackage{amscd} 
\usepackage{psfrag}
\usepackage{graphicx} 
\usepackage[latin1]{inputenc}  
\usepackage[all]{xy} 
\usepackage[mathcal]{eucal}

\theoremstyle{definition}

\theoremstyle{remark}

\def\interieur#1{\mathord{\mathop{\kern 0pt #1}\limits^\circ}}

\definecolor{NoteColor}{rgb}{1,0,0}

\title[On Grothendieck's tame topology]{On Grothendieck's tame topology}
\author[N. A'Campo, L. Ji and A. Papadopoulos]{Norbert A'Campo, Lizhen Ji \\
and Athanase Papadopoulos}

\address{N. A'Campo: Universit\"at Basel,  Mathematisches Institut, 
\\
Spiegelgasse 1, 4051 Basel, Switzerland
\\
 email:\,\tt{norbert.acampo@gmail.com}}
\address{L. Ji: Department of Mathematics, University of Michigan\\ Ann Arbor, MI 48109, USA
\\
email:\,\tt{lji@umich.edu}}
\address{A. Papadopoulos: Institut de Recherche Math\'ematique Avanc\'ee, UMR 7501
\\
Universit{\'e} de Strasbourg and CNRS,\\
7 rue Ren\'e Descartes, 67084 Strasbourg Cedex, France
\\
email:\,\tt{papadop@math.unistra.fr}}

 \date{\today}

\begin{document}

 \begin{abstract} 

Grothendieck's \emph{Esquisse d'un programme}\index{Grothendieck!Esquisse d'un programme}\index{Esquisse d'un programme (Grothendieck)}  is often referred to for the ideas it contains on dessins d'enfants, the Teichm\"uller tower, and the actions of the absolute Galois group on these objects or their \'etale fundamental groups; see the surveys \cite{AJP2} and \cite{Guillot} in the present volume. But this program contains several other important ideas. In particular, motivated by surface topology and moduli spaces of Riemann surfaces, Grothendieck calls there for a recasting of topology, in order to make it fit to the objects of semialgebraic and semianalytic geometry, and in particular to the study of the Mumford-Deligne compactifications of moduli spaces. A new conception of manifold, of submanifold and of maps between them is outlined. We review these ideas in the present chapter, because of their relation to the theory of moduli and Teichm\"uller spaces. We also mention briefly the relations between Grothendieck's ideas and  earlier theories developed by Whitney, \L{}ojasiewicz and Hironaka and especially Thom, and with the more recent  theory of o-minimal structures.

 \medskip


 \medskip


 \medskip

The final version of this paper will appear as a chapter in Volume VI of {\it the Handbook of Teichm\"uller theory}. This volume is dedicated to the memory of Alexander Grothendieck.

\end{abstract}

\bigskip

\maketitle

\tableofcontents

\section{Introduction} \label{s:intro}

Grothendieck's \emph{Esquisse d'un programme}\index{Grothendieck!Esquisse d'un programme}\index{Esquisse d'un programme (Grothendieck)} \cite{Gro-esquisse} is often referred to for the notions of dessin d'enfant and of Teichm\"uller tower,\index{Teichm\"uller tower}\index{tower!Teichm\"uller}\index{tower!Teichm\"ller} and for the actions of the absolute Galois group on these objects or their \'etale fundamental groups. But the \emph{Esquisse} also contains several other important ideas. In particular, motivated by surface topology and the theory of moduli spaces, Grothendieck calls there for a recasting of topology so that it becomes more adapted to the objects of semianalytic and semialgebraic geometry.
The name of the new field that he aims to is usually translated by \emph{tame topology} (cf. the English translation of \cite{Gro-esquisse}).\index{Grothendieck!Esquisse d'un programme}\index{Esquisse d'un programme (Grothendieck)} The French term that Grothendieck uses is ``topologie mod\'er\'ee" (moderate topology).\index{tame topology}\index{tame topology}\index{topology!tame}\index{topology!moderate} One of Grothendieck's aims is to obtain a branch of  topology which would give a  satisfactory theory of \emph{d\'evissage}\index{stratified structure!d\'evissage}\index{stratified structure!d\'evissage}\index{d\'evissage} (``unscrewing") of stratified structures. In the \emph{Esquisse}, he proposes some of the main axioms and he formulates the foundational  theorems of such a field. His motivation stems from the structure of the moduli spaces $\mathcal{M}_{g,n}$, a class of natural spaces which turn out to be subspaces of real analytic spaces, and maps between them. The moduli spaces are also the building blocks of the so-called modular tower,\index{tower!modular} which is a basic object in the \emph{Esquisse}; cf. \cite{AJP2}. The relation with Teichm\"uller spaces is the main reason why this chapter is included in the present Handbook.

           Stratified structures\index{stratified structure} (e.g. those underlying algebraic sets), rather than being given by atlases, are defined using attaching maps, and these attaching maps are required to have some regularity (for instance, Peano curves will not appear as attaching maps, the Jordan theorem should easily follow from the definitions, etc.).  Grothendieck proposed then a new branch of topology based on a system of axioms that would be a natural setting for semianalytic and semialgebraic geometries\index{geometry!semianalytic}\index{semianalytic!geometry}\index{geometry!semialgebraic}\index{semialgebraic!geometry} and that would rule out pathological situations like the one we mentioned.
           
  At the time Grothendieck wrote his \emph{Esquisse},\index{Grothendieck!Esquisse d'un programme}\index{Esquisse d'un programme (Grothendieck)} there was already a theory for stratified spaces; such a theory was first developed by Hassler Whitney,\index{Whitney, Hassler} see \cite{Whitney1} \cite{Whitney2}. It was generally admitted that this theory deals with the  stratification of an algebraic variety in an effective and algorithmic way. Ren\'e Thom\index{Thom, Ren\'e} also developed a theory of stratified spaces; cf. \cite{Thom2}, where he uses in particular his notion of \emph{controlled submersions}.  Thom's theory applies to finite complexes, manifolds with corners and semianalytic sets. However, it appears from the \emph{Esquisse} that Grothendieck considered that the existing theories were not flexible enough to include the study of objects like moduli spaces. For instance, in Thom's theory in \cite{Thom2}, a stratified set is obtained by gluing a \emph{finite} union of $C^\infty$ manifolds. This is obviously not the case for the augmented Teichm\"uller space boundary, which covers the one of the Mumford-Deligne compactification of moduli of curves, whose natural decomposition into a stratified space is not locally finite. In any case, a major issue in all these theories concerns the description of the attaching maps since these theories involve several complications regarding these maps. This is at the origin of Grothendieck's remark that classical topology, with the pathological cases that it necessarily involves, is not adapted to the examples he had in mind.

More recently, some new axiomatic topological theories were developed, which are adapted to the setting of semialgebraic geometry and also with the aim of ruling out the pathological  phenomena that one encounters in the setting of classical topology. One of these theories, which is considered to be in the lineage of Grothendieck's ideas, is  that of o-minimal structures.\index{o-minimal structure} This theory is based on a certain number of axioms which specify the kind of subsets of $\mathbb{R}^n$ that are accepted, and the functions are defined as those whose graphs belong to the admissible sets. Certain authors (cf. \cite{den-Dries} and \cite{D1}) consider that this is the theory which makes precise Grothendieck's ideas of tame topology. There is a theory, in mathematical logic, which carries the same name (o-minimal theory), and which is part of the general theory of \emph{quantifier elimination}\index{quantifier elimination} in model theory.\index{model theory}  In the the theory of quantifier elimination, one tries to replace (mathematical) formulae which contains quantifiers ($\forall, \exists$) by formulae without quantifiers and admitting the same models. One example, in the theory of real fields (that is, the fields that have the same first-order properties as the field of real numbers), is the formula $\exists x,  x^2+bx+c=0$ which is equivalent to the formula $b^2-ac\geq 0$ which uses no quantifiers.
The two theories (the geometric o-minimal theory and o-minimal theory as a branch of mathematical logic) are related. We shall say more about that in the last section of this chapter. Let us recall that the combination of geometry with mathematical logic is not a new subject; it can  be traced back (at least) to the works of Hilbert and Frege.

In the following, we shall review Grothendieck's considerations. 
They are mainly contained in \S\,5 and \S\,6 of the \emph{Esquisse}.\index{Grothendieck!Esquisse d'un programme}\index{Esquisse d'un programme (Grothendieck)} We shall mention relations with the  theory of stratified sets developed by Thom. Thom's theories already contain an important part of Grothendieck's program on this subject (before Grothendieck formulated it), but without the number-theoretic background and motivations.

\bigskip

\noindent{\bf Acknowledgements.---} We would like to thank Bernard Teissier who read several preliminary versions of this chapter and made corrections. Most of the remarks on the work of Thom are due to him. We also thank  Mark Goresky and Fran\c cois Laudenbach for correspondence on this subject.

\section{Grothendieck's ideas on tame topology in the \emph{Esquisse}}
Regarding\index{Grothendieck!Esquisse d'un programme}\index{Esquisse d'un programme (Grothendieck)} his program on tame topology, Grothendieck writes: ``I was first and foremost interested by the modular algebraic multiplicities, over the absolute base-field $\mathbb{Q}$,
 and by a `d\'evissage'\index{stratified structure!d\'evissage}\index{d\'evissage} at infinity of their geometric fundamental groups (i.e. of the profinite Teichm\"uller groups) which would be compatible with the natural operations of $\Gamma=\mathrm{Gal}(\overline{\mathbb{Q}}/\mathbb{Q})$."\footnote{All the quotes that we make here are from \S\,5 and \S\,6 of the\index{Grothendieck!Esquisse d'un programme}\index{Esquisse d'un programme (Grothendieck)} \emph{Esquisse}. The English translation is due to P. Lochak and L. Schneps.} He then comments on the difficulties that one encounters as soon as he wants to make big changes in what is commonly considered as basic mathematics: 
 \begin{quote}\small What is again lacking is not the technical virtuosity of mathematicians, which is sometimes impressive, but the audacity (or simply the innocence... ) to free oneself from a familiar context accepted by a flawless consensus...
 \end{quote}
 
Grothendieck recalls that the field of topology at the time he wrote his \emph{Esquisse} was still dominated by the development, done during the 1930s and 1940s, by analysts, in a way that fits their needs, rather than by geometers. He writes that the problem with such a development is that one has to deal with several pathological situations that have nothing to do with geometry. He declares that the fact that ``the foundations of topology are inadequate is manifest from the very beginning, in the form of `false problems'  (at least from the point of view of the  topological intuition of shape)." These false problems include the existence of wild phenomena (space-filling curves, etc.) that add complications which are not essential. He states that a new field of topology is needed, one which should be adapted to a theory of ``d\'evissage"\index{stratified structure!d\'evissage}\index{d\'evissage} (unscrewing) of stratified structures, a device which he was led to use several times in his previous works. Stratifications naturally appear in real or complex analytic geometry, where singular sets of maps appear as decreasing sequences of nested singular loci of decreasing dimensions. For Grothendieck, moduli spaces of geometric structures are naturally stratified sets, and the stratification also appears in the degeneration theory of these  structures. The main examples are the moduli spaces $\mathcal{M}_{g,n}$ of algebraic curves, equipped with their Mumford-Deligne boundaries.  Grothendieck calls these spaces the Mumford-Deligne multiplicities, and he denotes them by 
$\widehat{\mathcal{M}}_{g,\nu}$.

It is interesting to stop for a while on the word ``multiplicity,"\index{multiplicity} which is a term used by Riemann and which is an ancestor of our word ``manifold."\footnote{The German word is \emph{Mannigfaltigkeit};  it is the word still used today to denote a manifold, and it is close to the English word ``manifold." But the French term, which Grothendieck uses in the\index{Grothendieck!Esquisse d'un programme}\index{Esquisse d'un programme (Grothendieck)} \emph{Esquisse}, was replaced by the word ``vari\'et\'e," and is no more used to denote a mathematical object.} In fact, after Riemann, Poincar\'e and others used this word in its french form (``multiplicit\'e") to describe moduli space. Their view on that space was close to what we intend today by a manifold, given as a subset of a Euclidean space defined by a certain number of equations. Poincar\'e mentioned explicitly that the  moduli spaces of Riemann surfaces have singularities (although he did not specify the nature of these singularities). We refer the interested reader to the historical survey \cite{2015a} and to the more specialized one \cite{Oh}, for the notion of multiplicity and the origin of the notion of manifold.\index{multiplicity} At several places in the \emph{Esquisse},\index{Grothendieck!Esquisse d'un programme}\index{Esquisse d'un programme (Grothendieck)}  Grothendieck explains the meaning he gives to the  word ``multiplicity," with a reference to a stratification indexed by graphs that parametrize the possible combinatorial structures of stable curves (\S\,5). In fact, his use of the word ``multiplicity" is close to what we call today an orbifold: 
\begin{quote}\small
Two-dimensional geometry provides many other examples of such modular stratified structures, which all (if not using rigidification) appear as ``multiplicities" rather than as spaces or manifolds in the usual sense (as the points of these multiplicities may have non-trivial automorphism groups).
\end{quote}

 Other examples of stratified spaces that arise from the geometry of surfaces mentioned by Grothendieck are polygons (and he specifies, Euclidean, spherical or hyperbolic), systems of straight lines in a projective plane, systems of ``pseudo-straight lines" in a projective topological plane, and more general immersed curves with normal crossings. 
 
 Beyond these examples, Grothendieck declares he had ``the premonition of the ubiquity of stratified structures in practically all domains of geometry." One should note that a similar idea was expressed by Thom, who considered that all the sets that one encounters in geometry (in particular in generalized singular loci) are, at least in generic stable situations, stratified sets.
   
 As the simplest example of a stratified structure, Grothendieck mentions  pairs $(X,Y)$ of topological manifolds, where $X$ is a closed submanifold of $Y$ such that $X$ along $Y$ is ``equisingular."\index{equisingular set} In order to fit into this theory of stratification, such a notion needs to be defined with care, and for that purpose, one has to specify in a precise way the kind of \emph{tubular neighborhood}\index{tubular neighborhood} of $X$ in $Y$ that is needed in the presence of such a structure which is more rigid than the topological one. The classical cases are  the piecewise-linear, the Riemannian, and the metric. Such a tubular neighborhood has to be canonical, that is, well defined up to an automorphism of the structure. Grothendieck says that for that purpose, one has to work systematically in the \emph{isotopic categories}\index{isotopic category}\index{category!isotopic} associated with the categories of topological nature that arise a priori.  Isotopic categories are categories where two maps are considered to be the same if they are isotopic. Let us mention again Thom,\index{Thom, Ren\'e} who, in the early 1960s, developed a notion of tubular neighborhood in the setting of stratified spaces,\index{stratified space}\index{space!stratified} and a theory of locally trivial stratification\index{locally trivial stratification} \cite{Thom1}. In particular, his two theorems on isotopy provide a template for a deep reflection on stratified sets and morphisms. For an overview of the work of Thom on this subject, we refer the reader to the report \cite{TT} by Teissier.

As already said, Grothendieck calls a topology which avoids the  pathological situations a \emph{tame topology},\index{tame topology}\index{topology!tame} (in French, ``topologie mod\'er\'ee") and he declares that such a topology, which he wishes to develop, will not be unique, but that there is a ``vast infinity" of possibilities. They range ``from the strictest of all, the one which deals with the `piecewise-algebraic spaces' (with $\overline{\mathbb{Q}}_r=\overline{\mathbb{Q}}\cap \mathbb{R}$) to the piecewise-analytic."\index{topology!piecewise-analytic}\index{piecewise-analytic!topology} One must mention here the theory of subanalytic spaces\index{space!subanalytic}\index{subanalytic space} developed by \L{}ojasiewicz and Hironaka, and in particular the latter's work on ``resolution of singularities."  In this theory, objects are defined not only by analytic equations, but also by analytic inequalities (they may have ``corners" in the analytic sense an even more complicated singularities). The term ``subanalytic" was introduced by Hironaka in his paper \cite{H73}. The theory dealing with these objects is rich.  There is a so-called \emph{uniformization theorem}\index{uniformization theorem!subanalytic sets}\index{subanalytic set!uniformization theorem} for subanalytic sets, which is a consequence of Hironaka's theory of resolution of singularities \cite{H64} \cite{H74}. This theorem says that any closed subanalytic subset of $\mathbb{R}^n$ is the image by an analytic map of a proper analytic manifold of the same dimension.

There are also algebraic versions of this theory. An important observation in this setting is that projections of algebraic varieties (say, over the reals) onto affine subspaces are defined by inequations, and not only by equations. The standard example is that the projection of the circle in $\mathbb{R}^2$ defined by the equation $x^2+y^2-1=0$ on the $x$-axis is the interval $[0,1]$. This is not an algebraic set. With this in mind, a subset of $\mathbb{R}^n$ is called  \emph{semialgebraic}\index{semialgebraic set} if it can be obtained using combinations (finite operations of unions, intersections and complements) of polynomial equations and polynomial inequalities. (If one uses only polynomial equations -- with no inequalities -- then one gets the usual definition of an \emph{algebraic set}.\index{algebraic set})  A theorem of \L{}ojasiewicz  \cite{L1} says that a semialgebraic set can be triangulated, that is, transformed into a linearly embedded simplicial complex by a semialgebraic map of the ambient space. His main tool was an inequality known as the \L{}ojasiewicz inequality. It says that if $U$ is an open subset of a Euclidean space $\mathbb{R}^n$ and   $f:U\to\mathbb{R}$ an analytic function, then for any compact subset $K\subset U$ there exists $\alpha>0, C>0$ such that for any $p$ in $K$, $d(p,Z)^\alpha\leq C\vert f(p)\vert$, where $Z$ is the analytic subset of $U$ where $f$ vanishes. See also \cite{L2} for the work of  \L{}ojasiewicz.

 The Cartesian product of two semialgebraic sets is semialgebraic. Stratified sets naturally appear in this theory, since every semialgebraic set admits a stratification by semialgebraic sets of decreasing dimension. In particular, the boundary of a semialgebraic set is a semialgebraic set of lower dimension. These properties and others are summarized in the report \cite{Marker}.

The \emph{Tarski-Seidenberg theorem}\index{Tarski-Seidenberg theorem}\index{theorem!Tarski-Seidenberg} says that the projection onto an affine subspace of a semialgebraic set is semialgebraic. This theorem also implies that the closure of a semialgebraic set is semialgebraic. In fact, for any semialgebraic subset $X$ of $\mathbb{R}^n$, the closure $\overline{X}$
of $X$, its interior, and its frontier are semi-algebraic sets. (See \cite{Benedetti} Propostion 2.3.7.) The Tarski-Seidenberg theorem is also discussed in the books \cite{Khovanski} and \cite{Bochnak}, and it is summarized, with several other related things, in the review \cite{Marker}. From the point of view of model theory, this theorem is an illustration of the theory of quantifier elimination over the reals.

After talking about semialgebraic sets, one needs to talk about maps between them. There are some natural properties that such maps must satisfy, and this leads to the notion  of isomorphism between two semialgebraic sets. 

Algebraic and semialgebraic sets are tame in the sense that the pathological examples of Cantor sets, space filling curves, Sierpi\'nski sponges, etc. do not occur as level sets of algebraic or semialgebraic functions. In some sense, this justifies Grothendieck's assertions in the  \emph{Esquisse},\index{Grothendieck!Esquisse d'un programme}\index{Esquisse d'un programme (Grothendieck)}  that such pathologies are irrelevant in this kind of geometry.

In his program, 
Grothendieck states some of the foundational theorems he expects to hold in tame topology. One of them is a \emph{comparison theorem}, which says the following:
 \begin{quote}
 \emph{We essentially find the same isotopic categories (or even  $\infty$-isotopic) whatever the tame theory we work with.}
 \end{quote}
  He makes more precise statements about this idea in \S\,5 of the \emph{Esquisse}. The maps considered in this category may be embeddings, fibrations, smooth, \'etale fibrations, etc. Among the axioms that he introduces is the \emph{triangulability axiom}\index{triangulability axiom} ``in the tame sense, of a tame part of $\mathbb{R}^n$." He considers a  ``piecewise $\mathbb{R}$-algebraic" theory of complex algebraic varieties, and his setting also includes varieties defined over number fields. Again, one has to mention here that Thom considered that any semianalytic set, in a neighborhood of any of its points, is equivalent, by an ambient isotopy, to a semialgebraic set. Several variants of this fact were proved, in particular by T. Mostowski.

After the comparison theorem, the next fundamental theorem that Grothen\-dieck mentions concerns the existence and uniqueness of a tubular neighborhood $T$ for a closed tame subspace $Y$ in a tame space $X$, together with a way of constructing it, using for instance a tame map $X\to\mathbb{R}^+$ having $Y$ as a zero set, and the description of the boundary (not in the usual sense of a manifold with boundary) $\partial T$ of $T$. This is where the ``equisingularity" hypothesis on $X$ is needed.\index{equisingular set}  One expects  that the tubular neighborhood $T$ will be endowed, ``in an essentially unique way," with the structure of a locally trivial fibration over $Y$, with $\partial T$ as a subfibration. In fact, this is one of the isotopy theorems of Thom. In this respect, one of the basic ideas of Thom concerns tubular neighborhoods of strata; it says that each stratum has a ``tubular neighborhood" in the union of adjacent strata, and that the isotopies between these various neighborhoods are insured by appropriate transversality conditions.

Concerning his project of tame topology, Grothendieck writes: ``This is one of the least clear points in my temporary intuition of the situation, whereas the homotopy class of the predicted structure map $T\to Y$ has an obvious meaning, independent of any equisingularity hypothesis, as the homotopic inverse of the inclusion map $Y\to T$, which must be a homotopism."
He declares (\cite{Gro-esquisse} \S 5): 
\begin{quote}\small It will perhaps be said, not without reason, that all this may be only dreams, which will vanish in smoke as soon as one sets to work on specific examples, or even before, taking into account some known or obvious facts which have escaped me. Indeed, only working out specific examples
will make it possible to sift the right from the wrong and to reach the true substance. The only thing in all this which I have no doubt about, is the very necessity of such a foundational work, in other words, the artificiality
of the present foundations of topology, and the difficulties which they cause at each step. It may be however that the formulation I give of a theory
of d\'evissage\index{stratified structure!d\'evissage}\index{d\'evissage} of stratified structures in terms of an equivalence theorem of
suitable isotopic (or even $\infty$-isotopic) categories is actually too optimistic. But I should add that I have no real doubts about the fact that the theory
of these d\'evissages which I developed two years ago, although it remains in part heuristic, does indeed express some very tangible reality.
\end{quote}

For what concerns stratified spaces, the theory of tubular neighborhoods that Grothendieck sketches is included in a more general theory of ``local retraction data which  make it possible to
construct a canonical system of spaces, parametrized by the ordered set of
flags $\mathrm{Fl}(I)$ of the ordered set $I$ indexing the strata; these spaces [...] are connected by embedding and
proper fibration maps, which make it possible to reconstitute in an equally canonical way the original stratified structure, including these additional
structures." The main examples are again the Mumford-Deligne multiplicities\index{multiplicity!Mumford-Deligne} $\widehat{\mathcal{M}}_{g,\nu}$ with their canonical stratification at infinity.  Grothendieck writes  (\cite{Gro-esquisse} \S 5):
\begin{quote}\small
Another, probably less serious difficulty, is that this so-called moduli ``space" is in fact a multiplicity\index{multiplicity} -- which can be technically expressed by the necessity of replacing the index set $I$ for the strata with an (essentially finite) category of indices, here the ``MD [Mumford-Deligne] graphs" which ``parametrize" the possible ``combinatorial structures" of a stable curve of type $(g,\nu)$. This said, I can assert that the general theory of d\'evissage,\index{stratified structure!d\'evissage}\index{d\'evissage} which has been developed especially to meet the needs of this example, has indeed proved to be a precious guide, leading to a progressive understanding, with flawless coherence, of some essential aspects of the Teichm\"uller tower\index{tower!Teichm\"uller} (that is, essentially the ``structure at infinity" of the ordinary Teichm\"uller groups). It is this approach which finally led me, within some months, to the principle of a purely combinatorial construction of the tower of Teichm\"uller groupoids.\index{Teichm\"uller groupoid}
\end{quote}

In \S\,6, Grothendieck writes that one of the most interesting foundational theorems in that theory would be a d\'evissage theorem for maps $f:X\to Y$, where $Y$ is equipped with a filtration $Y^i$ by closed tame subspaces and where above the strata $Y^i\setminus Y^{i-1}$, $f$ induces a trivial fibration (from this tame point of view). This theorem may also be generalized to the case where the space $X$ is also equipped with a filtration.  He declares that theories of locally and globally tame spaces, and of set-theoretic differences of tame spaces, and of globally tame maps, generalizing the notion of locally trivial fibration, must be developed. Again, we note that part of this program was realized by Thom, cf. \cite{Thom2}.
 
%
\section{Grothendieck's later comments}
A few years after he wrote the \emph{Esquisse},\index{Grothendieck!Esquisse d'un programme}\index{Esquisse d'un programme (Grothendieck)} Grothendieck, commenting on his ideas on tame topology, considered that they did not attract the attention of topologists or geometers. In Section 2.12, Note 42 of his \emph{R\'ecoltes et semailles}\index{Grothendieck!R\'ecoltes et semailles}\index{R\'ecoltes et semailles (Grothendieck)} (1986) \cite{RS}, he compares the introduction of this new ``tame topology" to the introduction of schemes in the field of  algebraic geometry. He writes:\footnote{The translations from  \emph{R\'ecoltes et semailles} are ours.}
\begin{quote}\small A multitude of new invariants, whose nature is more subtle than the invariants which are currently known and used, but which I feel as fundamental, are planned in my program on ``moderate topology" (of which a very rough sketch is contained in my \emph{Sketch of a program}, which will appear in Volume 4 of the Reflections). This program is based on a notion of ``moderate topology," or of ``moderate space," which constitute, may be like the notion of topos, a (second) ``metamorphosis of the notion  of space."  It seems to me that it is much more evident and less profound than the latter. Nevertheless, I suspect that its immediate impact on topology itself will be much stronger, and that it will transform from top to bottom the craft of the topologist-geometer, by a profound transformation  of the conceptual context in which he works. (This also occurred in algebraic geometry, with the introduction of the point of view of schemes.) On the other hand, I sent my \emph{Esquisse}\index{Grothendieck!Esquisse d'un programme}\index{Esquisse d'un programme (Grothendieck)} to several old friends and eminent topologists, but it appears it did not have the gift of making anyone of them interested...\footnote{[Une foule de nouveaux invariants, de nature plus subtile que les invariants actuellement connus et utilis\'es, mais que je sens fondamentaux, sont pr\'evus dans mon programme de ``topologie mod\'er\'ee" (dont une esquisse tr\`es sommaire se trouve dans l'\emph{Esquisse d'un Programme}, \`a para\^\i tre dans le volume 4 des R\'eflexions). Ce programme est bas\'e sur la notion de "th\'eorie mod\'er\'ee" ou ``d'espace mod\'er\'e," qui constitue, un peu comme celle de topos, une (deuxi\`eme) ``m\'etamorphose de la notion d'espace." Elle est bien plus \'evidente (me semble-t-il) et moins profonde que cette derni\`ere. Je pr\'evois que ses retomb\'ees imm\'ediates sur la topologie ``proprement dite" vont \^etre pourtant nettement plus percutantes, et qu'elle va transformer de fond en comble le ``m\'etier" de topologue g\'eom\`etre, par une transformation profonde du contexte conceptuel dans lequel il travaille. (Comme cela a \'et\'e le cas aussi en g\'eom\'etrie alg\'ebrique avec l'introduction du point de vue des sch\'emas.) J'ai d'ailleurs envoy\'e mon ``Esquisse" \`a plusieurs de mes anciens amis et illustres topologues, mais il ne semble pas qu'elle ait eu le don d'en int\'eresser aucun...]}\end{quote}

In \S 15 Note $91_3$  of\index{Grothendieck!R\'ecoltes et semailles}\index{R\'ecoltes et semailles (Grothendieck)}  \emph{R\'ecoltes et semailles} \cite{RS}, he writes: 
\begin{quote}\small It is especially since my talks at the Cartan seminar on the foundations of complex analytic spaces and on the precise geometric interpretation of ``modular varieties with level" \`a la Teichm\"uller, around the end of the 1950s, that I understood  the importance of a double generalization of the usual notions of ``manifold" with which we had worked until now (algebraic, real or complex analytic, differentiable -- and later on, their `moderate topology' variants). The one consisted in enlarging the definition in such a way that arbitrary ``singularities" are admitted, as well as nilpotent elements in the structure sheaf of ``scalar functions" -- modelled on my foundational work with the notion of scheme.
The other extension is towards a ``relativisation" over appropriate locally ringed toposes (the ``absolute" notions being obtained by taking as a basis a point topos). This conceptual work, which is mature since more than twenty-five years and which started in the thesis of Monique Hakim, is still waiting to be resumed. A particularly interesting case is the  notion of relative rigid-analytic space which allows the consideration of ordinary complex analytic spaces and rigid-analytic spaces over local rings with variable residual characteristics, like the ``fibers" of a unique relative rigid-analytic space; in the same way as the notion of relative scheme (which was eventually generally accepted) allows to relate to each other algebraic varieties defined over fields of different characteristics.\footnote{[C'est surtout depuis mes expos\'es au S\'eminaire Cartan sur les fondements de la th\'eorie des espaces analytiques complexes, et sur l'interpr\'etation g\'eom\'etrique pr\'ecise des ``vari\'et\'es modulaires \`a niveau" \`a la Teichm\"uller, vers la fin des ann\'ees cinquante, que j'ai compris l'importance d'une double g\'en\'eralisation des notions courantes de ``vari\'et\'e" avec lesquelles on a  travaill\'e jusqu'\`a pr\'esent (alg\'ebrique, analytique r\'eelle oucomplexe, diff\'erentiable -- ou par la suite, leurs variantes en ``topologie mod\'er\'ee"). L'une consiste \`a \'elargir la d\'efinition de sorte \`a admettre des ``singularit\'es" arbitraires, et des \'el\'ements nilpotent dans le faisceau structural des ``fonctions scalaires" --  sur le mod\`ele de mon travail de fondements avec la notion de sch\'ema. L'autre extension est vers une ``relativisation" au-dessus de topos localement annel\'es convenables (les notions ``absolues" \'etant obtenues en prenant comme base un topos ponctuel). Ce travail conceptuel, m\^ur depuis plus de vingt-cinq ans et amorc\'e dans la th\`ese de Monique Hakim, attend toujours d'\^etre repris. Un cas particuli\`erement int\'eressant est celui d'une notion d'espace rigide-analytique relatif, qui permet de consid\'erer des espaces analytiques complexes ordinaires et des espaces rigide-analytiques sur des corps locaux \`a caract\'eristiques r\'esiduelles variables, comme les ``fibres" d'un m\^eme espace rigide-analytique relatif ; tout comme la notion de sch\'ema relatif (qui a fini par entrer dans les m\oe{}urs) permet de relier entre elles des vari\'et\'es alg\'ebriques d\'efinies sur des corps de caract\'eristiques diff\'erentes.]}
\end{quote}

In the same manuscript (\cite{RS}  \S\,18),\index{Grothendieck!R\'ecoltes et semailles}\index{R\'ecoltes et semailles (Grothendieck)} Grothendieck returns to his subject; he writes:  
\begin{quote}\small There is a fourth direction, carried on during my past as a mathematician, which is directed towards a renewal from top to bottom of an existing field. This the ``moderate topology" approach in topology, on which I somehow elaborated in the \emph{Esquisse d'un Programme} (Sections 5 and 6). Here, as many times since my very faraway high-school years, it seems that I am still the only one to realize the richness and the emergency of a work on foundations that has to be done, and whose need here seems to me however more evident than ever. I have the clear impression that the development of the point of view of moderate topology, in the spirit alluded to in the \emph{Esquisse d'un programme}, would represent for topology a renewal whose scope is comparable to the one the theory of schemes brought in algebraic geometry, and without requiring an energy investment of a comparable size. Furthermore, I think that such a moderate topology will end up being a valuable tool in the development of arithmetic geometry, in particular in order to be able to prove ``comparison theorems" between the ``profinite" homotopic structure associated to a stratified scheme of finite type over the field of complex numbers (or, more generally, to a scheme-stratified multiplicity of finite type over this field), and the corresponding  ``discrete" homotopic structure, defined using a transcendental way, and up to appropriate (in particular, equisingularity) hypotheses.\index{equisingular set} This question only makes sense in terms of a precise ``d\'evissage theory"\index{stratified structure!d\'evissage}\index{d\'evissage} for stratified structures which seems to me, in the case of ``transcendental" topology, to require the introduction of the ``moderate" context."\footnote{[Il y a enfin une quatri\`eme direction de r\'eflexion, poursuivie dans mon pass\'e de math\'ematicien, allant en direction d'un renouvellement ``de fond en comble" d'une discipline existante. Il s'agit de l'approche "topologie mod\'er\'ee" en topologie, sur laquelle je m'\'etends quelque peu dans l'``Esquisse d'un Programme" (par. 5 et 6). Ici, comme tant de fois depuis les ann\'ees lointaines du lyc\'ee, il semblerait que je sois seul encore \`a sentir la richesse et l'urgence d'un travail de fondements \`a faire, dont le besoin ici me para\^\i t plus \'evident pourtant que jamais. J'ai le sentiment tr\`es net que le d\'eveloppement du point de vue de la topologie mod\'er\'ee, dans l'esprit \'evoqu\'e dans l'Esquisse d'un programme, repr\'esenterait pour la topologie un renouvellement de port\'ee comparable \`a celui que le point de vue des sch\'emas a apport\'e en g\'eom\'etrie alg\'ebrique, et ceci, sans pour autant exiger des investissement d'\'energie de dimensions comparables. De plus, je pense qu'une telle topologie mod\'er\'ee finira par s'av\'erer un outil pr\'ecieux dans le d\'eveloppement de la g\'eom\'etrie arithm\'etique, pour arriver notamment \`a formuler et \`a prouver des ``th\'eor\`emes de comparaisons" entre la structure homotopique ``profinie" associ\'ee \`a un sch\'ema stratifi\'e de type fini sur le corps des complexes (ou plus g\'en\'eralement, \`a une multiplicit\'e sch\'ematique stratifi\'ee de type fini sur ce corps), et la structure homotopique ``discr\`ete" correspondante, d\'efinie par voie transcendante, et modulo des hypoth\`eses (d'\'equisingularit\'e notamment) convenables. Cette question n'a de sens qu'en termes d'une ``th\'eorie de d\'evissage" pr\'ecise pour les structures stratifi\'ees, qui dans le cadre de la topologie ``transcendante" me semble n\'ecessiter l'introduction du contexte ``mod\'er\'e."]}
\end{quote}

In the section called \emph{La vision -- ou douze th\`emes pour une harmonie} (``The vision -- or twelve themes for a harmony") of \emph{R\'ecoltes et semailles} (\cite{RS} \S\,2.8), the subject of tame topology is  considered by Grothendieck as one of the twelve themes which he describes as his ``great ideas" (\emph{grandes id\'ees}).

\section{O-minimal sets}
In this section, we briefly mention a few facts on o-minimal structures.

The notion of o-minimal structure\index{o-minimal structure} is a kind of a generalization of a semialgebraic and semianalytic structure, and it is in part motivated by it.\footnote{According to Teissier, the theory is also motivated by  the work of \L{}ojasiewicz on semianalytic sets. The very definition of o-minimality generalizes the crucial observation made by \L{}ojasiewicz that the proofs of (most of) the tameness properties do not really require the global finiteness of generic projections of an algebraic subset to an affine space of the same dimension (Noether's theorem) that appears in the Tarski-Seidenberg theory, but only the local finiteness expressed by the Weierstrass preparation theorem.} We already mentioned that from their very definition, semialgebraic sets are stable under the usual Boolean operations of intersection, union and taking the complement. More precisely, an o-minimal structure\index{o-minimal structure} on $\mathbb{R}$ is a collection of subsets $S_n$ of $\mathbb{R}^n$ for each $n\geq 1$ satisfying the following:
\begin{itemize}
\item Each $S_n$  is stable under the operations of finite union, intersection, and taking the complement;
\item the elements of the collection $S_1$ are finite unions of intervals and points;
\item the projection maps from $\mathbb{R}^{n+1}$ to $\mathbb{R}^n$ sends subsets in $S_{n+1}$ to subsets in $S_n$. 
\end{itemize}
We already recalled that the Tarski-Seidenberg theorem\index{theorem!Tarski-Seidenberg} says that the projection to a lower-dimensional affine space of a semiagebraic set is semialgebraic, so that semialgebraic sets constitute an o-minimal structure. The study of o-minimal structures is also  a subfield of mathematical logic, and the theory can be wholly developed as a theory about quantifiers. It is closely related to model theory, even though the motivation behind it comes from the theory of semialgebraic sets. In fact, the theory of semialgebraic and subanalytic sets are prominent instances where the properties of quantifier elimination\index{quantifier elimination} in model theory\index{model theory} may be applied.

            Loi reports in \cite{Loi} that  the name o-minimal structure was given by van den Dries, Knight, Pillay and Steinhorn, who developed the general theory (\cite{D1}, \cite{KPS},  \cite{PS}). The author of \cite{Loi} also says that Shiota had a similar program (\cite{S1}, \cite{S2}). He then writes:   ``The theory of o-minimal structures is a wide-ranging generalization of semialgebraic and subanalytic geometry. Moreover, one can view the subject as a realization of GrothendieckÕs idea of topologie mod\'er\'ee, or tame topology, in his \emph{Esquisse dÕun Programme} (1984)."   The paper  \cite{S3} by Shiota is an example of the combination of topology, geometry and logic which is realized in o-minimal theory.

\section{As a way of conclusion}
We saw that some of the basic ideas of Grothendieck on tame topology were worked out by Thom and others, even before Grothendieck formulated his program. We also saw that there are relations of these ideas with other subjects that grew up after Grothendieck's work, like the theory of o-minimal structures. But Grothendieck's project of recasting the whole foundations of topology or of creating a new field of topology based on these ideas has still not been realized. We consider this project as another aspect of his broad vision on mathematics.

\bigskip

\noindent \emph{Acknowledgements.} The authors acknowledge support from U.S. National Science Foundation grants DMS 1107452, 1107263, 1107367 ``RNMS: GEometric structures And Representation varieties" (the GEAR Network). The authors wrote this chapter during a stay at  the Mittag-Leffler Institute, and they would like to thank the institute for its hospitality.

         \printindex

\begin{thebibliography}{XXXX}
 
 

\bibitem{2015a} N. A'Campo, L. Ji and A. Papadopoulos,  On the early history of moduli and Teichm\"uller spaces, In \emph{Lipman Bers, a life in mathematics}, a volume dedicated to Lipman Bers, ed. L. Keen, I. Kra and Rubi Rodriguez, American Math. Society, 2015, 175-262.
 
 \bibitem{AJP1} N. A'Campo, L. Ji and A. Papadopoulos, On Grothendieck's construction of Teichm\"uller space. In \emph{Handbook of Teichm\"uller theory}  (A. Papadopoulos, ed.), Volume VI, EMS Publishing House,  Z\"urich, 2016, 35-69.

 
  \bibitem{AJP2} N. A'Campo, L. Ji and A. Papadopoulos, Actions of the absolute Galois group. In \emph{Handbook of Teichm\"uller theory}  (A. Papadopoulos, ed.), Volume VI, EMS Publishing House,  Z\"urich, 2016, 397-435.

 
\bibitem{Benedetti} R. Benedetti and J.-J. Risler, Real algebraic and semi-algebraic sets, Hermann, 1990.

\bibitem{Bierstone} E. Bierstone and P. Milman, Subanalytic Geometry, In: \emph{Model Theory, Algebra, and Geometry},  (D.Haskell, A. Pillay, and Ch. Steinhorn, ed.) MSRI Publications, Volume 39, 2000

\bibitem{Bochnak} J. Bochnak, M. Coste and M.-F. Roy, \emph{Real Algebraic Geometry}, Springer Verlag, 1998.


\bibitem{den-Dries} L. van den Dries,  \emph{Tame topology and o-minimal structures}. London Mathematical Society lecture note series, no. 248. Cambridge University Press, Cambridge, New York, and Oakleigh, Victoria, 1998.

           \bibitem{D1}  L. van den Dries, A generalization of the Tarski-Seidenberg theorem and some nondefinability results, \emph{Bull. Amer. Math. Soc}. (N.S.) 15 (1986), 189-193. 
%
           
          
 \bibitem{Gro-esquisse} A. Grothendieck, Esquisse d'un programme (Sketch of a program), unpublished manuscript (1984), English translation by P. Lochak and  L. Schneps  in  \emph{Geometric Galois actions, vol. 1, ``Around Grothendieck's Esquisse d'un Programme"} (L. Schneps and P. Lochak, ed.) London Math. Soc. Lecture Note Ser. vol.  242, pp. 5-48, Cambridge Univ. Press, Cambridge, 1997. 
 
 
 \bibitem{RS} A. Grothendieck,  R\'ecoles et semailles : R\'eflexions et t\'emoignage sur un pass\'e de math\'ematicien (Harvesting and Sowing : Reflections and testimony on a mathematician's past), manuscript, 1983-1986, to appear, 



  \bibitem{Guillot} P. Guillot, A primer on dessins,   In \emph{Handbook of Teichm\"uller theory}  (A. Papadopoulos, ed.), Volume VI, EMS Publishing House,  Z\"urich, 2016.
  
  \bibitem{H64} H. Hironaka, Resolution of singularities of an algebraic variety over a field of characteristic zero, \emph{Ann. of Math}. (2) 79 (1964), 109-203, 205-326.
  
  
  \bibitem{H73} H. Hironaka, Subanalytic sets, In: \emph{Number theory, algebraic geometry and commutative algebra: in honor of Yasuo Akizuki}, (Y. Kusunoki  et al. ed.), Kinokuniya, Tokyo, 1973, p. 453-493.
  
    \bibitem{H74} H. Hironaka, \emph{Introduction to the theory of infinitely near singular points}, Mem. Mat. Instituto Jorge Juan 28, Consejo Superior de Investigaciones Cient\'\i ficas, Madrid, 1974.
  
  
  
 \bibitem{Khovanski} A. G. Khovanskii,  \emph{Fewnomials}. Translations of Mathematical Monographs 88. Translated from the Russian by Smilka Zdravkovska. Providence, RI: American Mathematical Society, 1991.
 
          \bibitem{KPS}  J. Knight, A. Pillay and C. Steinhorn, Definable sets in ordered structures II, \emph{Trans. Amer. Math. Soc.} 295 (1986), 593-605.
          

\bibitem{Loi-Phan} T. L. Loi and  P. Phan, \emph{Acta Mathematica Vietnamica}, 29 (2014)  No. 4,  637-647
 
\bibitem{Loi} T. L. Loi, o-minimal structures. In: \emph{The Japanese-Australian workshop on real and complex singularities} (Fukui, Toshizumi et al., ed.), The University of Sydney, Sydney, Australia, September 15-18, 2009.  Proceedings of the Centre for Mathematics and its Applications, Australian National University 43,  2010, 19-30.  


\bibitem{Loi2} T. L. Loi, Tame topology and Tarski-type systems, \emph{Vietnam J. Math.} 31 (2003), No. 2, 127-136.
 
            
 \bibitem{L1} S. \L{}ojasiewicz, Triangulation of semi-analytic sets, \emph{Ann. Scuola Norm. Sup. Pisa} (3), 18 (1964), 449-474. 


 \bibitem{L2} S. \L{}ojasiewicz, Sur les ensembles semi-analytiques. In \emph{Actes du Congr\`es International des Math\'ematiciens (Nice, 1970)}, Volume 2,  Gauthier-Villars, Paris, 1971, 237-241.



\bibitem{Marker} D. Marker, Book review of \emph{Tame topology and o-minimal structures} by Lou van den Dries, Cambridge Univ. Press, New York, 1998. \emph{Bull. Amer. Math. Soc.} (N.S.)  37 (2000) No. 3, 351-357.


 \bibitem{Oh} K. Ohshika, The origin of the notion of manifold,  In \emph{From Riemann to differential geometry and relativity}, (L. Ji, A. Papadopoulos and S. Yamada, ed.), Springer Verlag, to appear in 2017.
 
 
 \bibitem{PS} A. Pillay and C. Steinhorn, Definable sets in ordered structures I, \emph{Trans. Amer. Math. Soc.}, 295 (1986), 565-592.
           


 
 \bibitem{S1} M. Shiota, Geometry of subanalytic and semianalytic sets: Abstract, In: \emph{Real analytic and Agebraic Geometry}, Proceedings of the international conference, Trento, Italy, September 21-25, 1992. Walter de Gruyter, Berlin, 1995,  251-275.
    
 \bibitem{S2} M. Shiota, \emph{Geometry of Subanalytic and Semialgebraic Sets}, Progress in Math., vol. 150, Birkh\"auser, Boston, 1997.
 
 \bibitem{S3}  M. Shiota, O-minimal Hauptvermutung for polyhedra I. \emph{Invent. math.}  196 (2014), Issue 1, 163-232.
           
\bibitem{TT} B. Teissier, Travaux de Thom sur les singularit\'es, \emph{Publ. Math. Inst. Hautes \'Etudes Sci.} 68 (1988), 19-25.

\bibitem{Teissier} B. Teissier, Tame and stratified objects. In: \emph{Geometric Galois Actions, 1. Around GrothendieckÕs esquisse dÕun programme}. London. Math. Soc. Lecture Note Series 242, 1997, 231-242.

\bibitem{Thom1} R. Thom, La stabili\'e topologique des applications polynomiales, \emph{LÕEns. Math}. (2) 8 (1962), 24-33. 

\bibitem{Thom2} R. Thom, 
Ensembles et morphismes stratifi\'es.  
  \emph{Bull. Am. Math. Soc.} 75 (1969), 240-284.



\bibitem{Whitney1} H. Whitney, Complexes of manifolds, \emph{Proc. Nat. Acad. Sci. U.S.A.} 33, (1947), 10-11.
 
\bibitem{Whitney2} H. Whitney, Elementary structure of real algebraic varieties, \emph{Ann. Math.} 66 (1957), 545-556.

   

\end{thebibliography}
  \end{document}